\documentclass[11pt]{article}
\usepackage{bbm}
\usepackage{amsfonts}
\usepackage{graphicx}
\usepackage{mathptmx}
\usepackage{hyperref}
\usepackage{latexsym,amsmath,amssymb,amsfonts,amsthm}

\title{THE TOPOLOGICAL RIGIDITY THEOREM FOR SUBMANIFOLDS IN SPACE FORMS\footnote{Research supported by the National Natural Science Foundation of China, Grant Nos. 11531012, 11371315, 11771394.}}

\author{JUAN-RU GU AND HONG-WEI XU\footnote{Corresponding author}}
\date{}
\numberwithin{equation}{section} \textheight 220 mm \textwidth 150
mm \evensidemargin  4mm \oddsidemargin  4mm \topmargin -10mm

\begin{document}
\maketitle
\begin{abstract}
Let $M$ be an $n(\geq 4)$-dimensional compact submanifold
 in the simply connected space form
$F^{n+p}(c)$ with constant curvature $c\geq 0$, where $H$ is the mean curvature of $M$.
We verify that if the scalar curvature of $M$ satisfies
 $R>n(n-2)(c+H^2)$, and if $Ric_M\geq (n-2-\frac{2\sigma_n}{2n-\sigma_n})(c+H^2)$,
 then $M$ is homeomorphic to a sphere. Here $\sigma_n=sgn(n-4)((-1)^n+3)$, and $sgn(\cdot)$ is the standard sign function. This improves our previous sphere theorem \cite{XG2}. It should be emphasized that our pinching
conditions above are optimal. We also obtain some new topological sphere
theorems for submanifolds with pinched scalar curvature and Ricci curvature. 
\end{abstract}

 \section{Introduction}
 It plays an important role in global differential
geometry to investigate geometrical and topological structures of
manifolds. After the pioneering rigidity theorem for closed minimal
submanifolds in a sphere due to Simons \cite{Simons}, several
striking rigidity results for minimal submanifolds were proved by
Chern, do Carmo, Kobayashi, Lawson, Yau and others (see \cite{Chern,
Ejiri, Lawson, Li, Yau}, etc.). In 1979, Ejiri \cite{Ejiri} obtained
the following rigidity theorem for compact minimal submanifold in a unit sphere.\\\\
\textbf{Theorem A.} \emph{Let $M$ be an $n(\geq4)$-dimensional
  oriented compact simply connected minimal submanifold in $S^{n+p}$.
If the Ricci curvature of $M$ satisfies $Ric_M \geq n-2,$ then $M$
is either the totally geodesic submanifold $S^n$, the Clifford torus
$S^{m}\big(\sqrt{\frac{1}{2}}\big)\times
S^{m}\big(\sqrt{\frac{1}{2}}\big)$ in $S^{n+1}$ with $n=2m$, or
$\mathbb{C}P^{2}(\frac{4}{3})$ in $S^7$. Here
$\mathbb{C}P^{2}(\frac{4}{3})$ denotes the $2$-dimensional complex
projective space minimally immersed into $S^7$ with constant
holomorphic sectional curvature $\frac{4}{3}$.}\\

The pinching constant above is the best possible in even dimensions.
In 1987, Sun \cite{Sun} studied the rigidity of compact oriented
submanifolds of dimension $n\geq 4$ with parallel mean curvature and
pinched Ricci curvature in a sphere, and partially generalized
Ejiri's rigidity theorem. In 1990s, Shen \cite{ShYB} and Li \cite{Li2} extended Ejiri's
rigidity theorem to the case of 3-dimensional compact minimal
submanifolds in a sphere. Moverover,  Li
\cite {Li2} discussed $n(\geq5)$-dimensional
compact minimal submanifolds in $S^{n+p}$, and obtained that if $n$ is odd, and
$Ric_{M}>n-2-\frac{1}{n-1}$, then $M$ is the totally geodesic. In 2011, Xu and Tian \cite{XT} obtained a
refined version of the Ejiri rigidity theorem without the assumption
that $M$ is simply connected.

 Set $UM=\bigcup_{x\in M}U_xM$, where $U_xM=\{u\in T_xM:\|u\|=1\}$, and $T_xM$ is the tangent space at $x\in M$. In 1986, Gauchman \cite{Gauchman} proved the following rigidity theorem for compact minimal submanifolds in a unit sphere.\\\\
 \textbf{Theorem B.} \emph{Let $M$ be an n-dimensional compact minimal submanifold  in $S^{n+p}$. If $$||h(u,u)||^{2}\leq C(n,p)$$ for all $u\in UM$, where $h$ is the second fundamental form of $M$, then either $||h(u,u)||^{2}=0$, and $M$
is totally geodesic, or $n$ is even, and $||h(u,u)||^{2}=C(n,p)$.
 Here $$C(n,p)=\left\{\begin{array}{ll} 1,&\mbox{\ for\ } p=1,\mbox{\ and\ } n\mbox{\  is even},\\
\frac{n}{n-1}, &\mbox{\ for\ } p=1,\mbox{\ and\ } n\mbox{\  is odd},\\
 \frac{1}{3}, &\mbox{\ for\ }  p\geq2, \mbox{\ and\ } n\mbox{\  is even},\\
 \frac{n}{3n-2}, &\mbox{\ for\ }  p\geq2, \mbox{\ and\ } n\mbox{\  is odd}.
 \end{array}\right.$$}

When $n$ is even, the above rigidity theorem is optimal. In 1991, Leung \cite{Leung}  proved that if  $M^n$ is an odd-dimensional compact minimal submanifold  in $S^{n+p}$ with flat normal connection, and if
$||h(u,u)||^{2}\leq \frac{n}{n-1}$ for all $u\in UM$, then $M$ is totally geodesic. Let $M=S^{m}\big(\sqrt{\frac{m}{n}}\big)\times
S^{m+1}\big(\sqrt{\frac{m+1}{n}}\big)$ be a Clifford minimal hypersurface
in $S^{n+1}$ with $n=2m+1$, then $\frac{n-1}{n+1}\leq||h(u,u)||^{2}\leq \frac{n+1}{n-1}$. Based on these facts, Leung proposed the following conjecture.\\\\
\textbf{Conjecture C.} \emph{Let  $M^n$ be an odd-dimensional compact minimal submanifold  in $S^{n+p}$. \\
(i) If $||h(u,u)||^{2}\leq \frac{n}{n-1}$ for all $u\in UM$, then $M$ is homeomorphic to a sphere.\\
(ii) If $||h(u,u)||^{2}<\frac{n+1}{n-1}$ for all $u\in UM$, then $M$ is homeomorphic to a sphere.}\\

In 2001, Hasanis and Vlachos \cite{HV} proved if $M^n(n\geq5)$ is an odd-dimensional oriented compact minimal submanifold in $S^{n+p}$, and if $Ric_M>n-2-\frac{2}{n-1}$, then $M$ is homeomorphic to a sphere. Moreover, they showed that the condition $Ric_M>n-2-\frac{2}{n-1}$ is equivalent to $||h(u,u)||^{2}<\frac{n+1}{n-1}$ when $M$ is a compact minimal submanifold in $S^{n+p}$ with flat normal connection. In 2016, Qian and Tang \cite{QT} gave some counter examples  to Conjecture C, which are focal submanifolds of isoparametric hypersurfaces in unit spheres with $g=4$.

 Since 1970s, the sphere theorems for general submanifolds have been
investigated by several authors \cite{Lawson2,Shiohama1,Shiohama2,Vlachos,Xia,Xin,XHZ,XZ}. Recently, the authors \cite{XG2} generalized the Ejiri rigidity theorem for
minimal submanifolds in a sphere due to Ejiri \cite{Ejiri} and Xu-Tian \cite{XT} to compact
submanifolds with parallel mean curvature in space forms, and proved the following topological sphere theorem for compact submanifolds in a space form with nonnegative constant curvature.\\\\
\textbf{ Theorem D.} \emph{Let $M$ be an $n(\geq4)$-dimensional compact submanifold in
$F^{n+p}(c)$ with  $c\geq 0$,
 if  $Ric_M>(n-2)(c+H^2),$
then $M$ is homeomorphic to a sphere.} \\

Moreover, Xu-Leng-Gu \cite{XG3}  investigated the compact submanifolds with odd dimension in
  space forms, and proved the following topological
   sphere theorem. \\\\
\textbf{Theorem E.} \emph{Let $M$ be an $n(\geq5)$-dimensional
compact submanifold in $F^{n+p}(c)$ with  $c\geq0$. If $n$ is odd,
and
$$Ric_M>(n-2-\epsilon_n)(c+H^2),$$  then M
is homeomorphic to a sphere. Here $\epsilon_n =\frac{4}{n^3-2n^2-n-2}.$}\\

In this paper, we study the compact submanifolds in space forms, and prove the following topological sphere theorems.\\\\
\textbf{Theorem 1.1.} \emph{Let $M$ be an $n(\geq4)$-dimensional compact submanifold in
$F^{n+p}(c)$ with  $c\geq0$. If the scalar curvature of M satisfies  $R>n(n-2)(c+H^2)$, and $$Ric_M\geq (n-2-\frac{2\sigma_n}{2n-\sigma_n})(c+H^2),$$
then $M$ is homeomorphic to a sphere. Here $\sigma_n=sgn(n-4)((-1)^n+3)$, and $sgn(\cdot)$ is the standard sign function.}\\\\
\textbf{Remark.} If $Ric_M>(n-2)(c+H^2)$, then $R>n(n-2)(c+H^2)$. Therefore, Theorem 1.1 is a generalization of Theorem D.\\\\
  \textbf{Theorem 1.2.} \emph{Let $M$ be an $n(\geq4)$-dimensional compact oriented submanifold in
$F^{n+p}(c)$ with  $c\geq0$. If  $R>n(n-2)(c+H^2)$, and  $$Ric_M\leq (n-2+\frac{2\tau_n}{2n+\tau_n})(c+H^2),$$
then $M$ is homeomorphic to a sphere. Here $\tau_n=(-1)^n+3$.}\\

The following examples show that the pinching conditions in Theorems
1.1 and 1.2 are the best possible.\\\\
\textbf{Example 1.3.} Let
$S^{n+1}(\frac{1}{\sqrt{c+H^2}})$ be the totally umbilic sphere in
$F^{n+p}(c)$. Here $H$ is a nonnegative constant.

(i) Let
$M=S^{m}\big(\frac{1}{\sqrt{2(c+H^2)}}\big)\times
S^{m}\big(\frac{1}{\sqrt{2(c+H^2)}}\big)$ be a Clifford hypersurface
in $S^{n+1}(\frac{1}{\sqrt{c+H^2}})$ with $n=2m$ and $c+H^2>0$.  Then $M$ is a
compact submanifold in $F^{n+p}(c)$ with constant mean curvature $H$
and constant Ricci curvature $Ric_M\equiv (n-2)(c+H^2)$. Moreover, $R=n(n-2)(c+H^2).$

(ii) Let
$M=S^{m-1}\big(\sqrt{\frac{m-1}{n(c+H^2)}}\big)\times
S^{m+1}\big(\sqrt{\frac{m+1}{n(c+H^2)}}\big)$ be a Clifford hypersurface
in $S^{n+1}(\frac{1}{\sqrt{c+H^2}})$ with  $n=2m$ and  $c+H^2>0$.  Then $M$ is a
compact submanifold in $F^{n+p}(c)$ with constant mean curvature $H$,  $Ric(e_i)\equiv (n-2-\frac{4}{n-2})(c+H^2)$ for $1\leq i\leq m-1$, and $Ric(e_j)\equiv (n-2+\frac{4}{n+2})(c+H^2)$ for $m\leq j\leq n$. Moreover, $R=n(n-2)(c+H^2).$

(iii) When $n=2m+1$, let
$M=S^{m}\big(\sqrt{\frac{m}{n(c+H^2)}}\big)\times
S^{m+1}\big(\sqrt{\frac{m+1}{n(c+H^2)}}\big)$ be a Clifford hypersurface
in $S^{n+1}(\frac{1}{\sqrt{c+H^2}})$ with  $c+H^2>0$. Then $M$ is a
compact submanifold in $F^{n+p}(c)$ with constant mean curvature $H$,  $Ric(e_i)\equiv (n-2-\frac{2}{n-1})(c+H^2)$ for $1\leq i\leq m$, and $Ric(e_j)\equiv (n-2+\frac{2}{n+1})(c+H^2)$ for $m+1\leq j\leq n$. Moreover, $R=n(n-2)(c+H^2).$\\

\section{Notation and lemmas}
 Throughout this paper,
let $M$ be an $n$-dimensional compact Riemannian manifold
isometrically immersed into an $(n+p)$-dimensional complete and
simply connected space form $F^{n+p}(c)$. We shall make use of the
following convention on the range of indices:
$$ 1\leq A,B,C,\ldots\leq n+p;\ 1\leq i,j,k,\ldots\leq n;\ n+1\leq
\alpha,\beta,\gamma,\ldots\leq n+p.$$ We let \{$e_{A}$\} be local
orthonormal frames in $F^{n+p}(c)$ such that, restricted to $M$, the
$e_{i}$'s are tangent to \emph{M}. Let \{$\omega _{A}$\} and
\{$\omega _{AB}$\} be the dual frame field and the connection
1-forms of $F^{n+p}(c)$ respectively. Restricting these forms to
\emph{M}, we have
\begin{eqnarray}&&\omega_{\alpha i}=\sum_{j} h^{\alpha}_{ij}\omega_{j},  \, \,
h^{\alpha}_{ij}=h^{\alpha}_{ji},\nonumber\\
&&h=\sum_{\alpha,i,j}h^{\alpha}_{ij}\omega_{i}\otimes\omega_{j}\otimes
e_{\alpha},\,\,\xi=\frac{1}{n}\sum_{\alpha,i}h^{\alpha}_{ii}e_{\alpha},\nonumber
\\
&&R_{ijkl}=c(\delta_{ik}\delta_{jl}-\delta_{il}\delta_{jk})+\sum_{\alpha}(h^{\alpha}_{ik}h^{\alpha}_{jl}-h^{\alpha}_{il}h^{\alpha}_{jk}), \label{Gauss}\\
&&R_{\alpha\beta kl}=\sum_{i}(h^{\alpha}_{ik}h^{\beta}_{il}-h^{\alpha}_{il}h^{\beta}_{ik}),\label{Codazzi}
\end{eqnarray}
where $h, \xi, R_{ijkl},$ and $R_{\alpha\beta kl}$ are the second
fundamental form, the mean curvature vector, the curvature tensor
and the normal curvature tensor of $M$, respectively.  We define
$$S=|h|^{2}, \ H=|\xi|, \ H_{\alpha}=(h^{\alpha}_{ij})_{n\times
n}.$$
Denote by $Ric(u)$ the Ricci curvature of $M$ in direction of $u\in
UM$. From the Gauss equation, we have
\begin{equation}Ric(e_i)=(n-1)c
+\sum_{\alpha,j}\big[h_{ii}^{\alpha}h_{jj}^{\alpha}
                          -(h_{ij}^{\alpha})^2\big].\label{Ricci}\end{equation}
Set $Ric_{\min}(x)=\min_{u\in U_{x}M}Ric(u)$.
The scalar curvature $R$ of $M$ is given by
\begin{equation}R=n(n-1)c+n^{2}H^{2}-S.\label{Gauss2}\end{equation}
Choose an
orthonormal basis $\{e_{i}\}$ in $T_{x}M$ such that $e_{j_0}=u,\,\,
span\{e_{j_1},\ldots,e_{j_k}\}=V_{x}^k$, where the indices $1\le
j_0,j_1,\ldots,j_k\le n$ are distinct with each other.
 Using the notations of \cite{GX}, we set
\begin{eqnarray}&& Ric^{[s]}([e_{j_0},\ldots,e_{j_{n-1}}])  = \sum_ {p=0}^{s-1}\sum_
 {q=0}^{n-1}R_{j_pj_{q}j_pj_{q}},\nonumber\\
 &&Ric^{[s]}_{\min}(x)=\min_{
 \{e_{j_0},\ldots,e_{j_{n-1}}\}\subset T_{x}M} Ric^{[s]}([e_{j_0},\ldots,e_{j_{n-1}}]),\nonumber\\
 &&Ric^{[s]}_{\max}(x)=\max_{
 \{e_{j_0},\ldots,e_{j_{n-1}}\}\subset T_{x}M}
 Ric^{[s]}([e_{j_0},\ldots,e_{j_{n-1}}]).\label{2-5}\end{eqnarray}
  \textbf{Definition 2.1(\cite{GX}).} \emph{We call  $Ric^{[s]}([e_{j_0},\ldots,e_{j_{n-1}}])$ the $s$-th weak Ricci curvature of $M$, respectively.}\\

The nonexistence theorem for stable currents
in a compact Riemannian manifold $M$ isometrically immersed into
$F^{n+p}(c)$ is employed to eliminate the homology groups
$H_{q}(M;\mathbb{Z})$ for $0<q<n$, which was initiated by
Lawson-Simons \cite{Lawson2} and extended by Xin \cite{Xin}.\\\\
\textbf{Theorem 2.1.} \emph{Let $M^{n}$ be a compact submanifold in
$F^{n+p}(c)$ with $c\geq0$. Assume that$$
\sum_{k=q+1}^{n}\sum_{i=1}^{q}[2|h(e_{i},e_{k})|^{2}-\langle
h(e_{i},e_{i}),h(e_{k},e_{k})\rangle]<q(n-q)c$$ holds for any
orthonormal basis $\{e_{i}\}$ of $T_xM$ at any point $x\in M$, where
q is an integer satisfying $0<q<n$. Then there do not exist any
stable q-currents. Moreover,
$H_{q}(M;\mathbb{Z})=H_{n-q}(M;\mathbb{Z})=0,$ and $\pi_{1}(M)=0$
when $q=1$. Here $H_{i}(M;\mathbb{Z})$ is the $i$-th homology group
of M with integer
coefficients.}\\

To prove the sphere theorems for submanifolds, we need
to eliminate the fundamental group
$\pi_1(M)$ under the s-th weak Ricci curvature and the scalar curvature pinching condition, and get the following lemma.\\\\
  \textbf{Lemma 2.2.} \emph{Let $M$ be an $n(\geq4)$-dimensional
 compact submanifold in $F^{n+p}(c)$ with  $c\geq0$.
 Assume that one of the following  conditions holds: \\
  $(i)$ $R+2Ric_{\min}>n(n-1)(c+H^2),$\\
  $(ii)$  $3R-2Ric^{[n-1]}_{\max}>n(n-1)(c+H^2),$\\then $H_{1}(M;\mathbb{Z})=H_{n-1}(M;\mathbb{Z})=0,$ and $\pi_1(M)=0$.}\\\\
\textbf{Proof.} From (\ref{Gauss2}), we have
$$S-nH^2=n(n-1)(c+H^2)-R.$$ This together with (\ref{Ricci})
implies that
\begin{eqnarray} &&\sum_{k=2}^{n}[2|h(e_{1},e_{k})|^{2}-\langle
h(e_{1},e_{1}),h(e_{k},e_{k})\rangle]\nonumber\\
&=&2\sum_{\alpha}\sum_{k=2}^{n}(h^{\alpha}_{1k})^{2}-
\sum_{\alpha}\sum_{k=2}^{n}h^{\alpha}_{11}h^{\alpha}_{kk}\nonumber\\
&=&\sum_{\alpha}\sum_{k=2}^{n}(h^{\alpha}_{1k})^{2}-Ric(e_1)+(n-1)c\nonumber\\
&\leq&\frac{1}{2}(S-nH^2)-Ric(e_1)+(n-1)c\nonumber\\
&=&\frac{n(n-1)}{2}(c+H^2)-\frac{R}{2}-Ric(e_1)+(n-1)c.
\label{2-6} \end{eqnarray}
Since $$2Ric_{\min}+R>n(n-1)(c+H^2),$$ we have $$\sum_{k=2}^{n}[2|h(e_{1},e_{k})|^{2}-\langle
h(e_{1},e_{1}),h(e_{k},e_{k})\rangle]<(n-1)c.$$
On the other hand, we also get from (\ref{2-6}) that
\begin{eqnarray} &&\sum_{k=2}^{n}[2|h(e_{1},e_{k})|^{2}-\langle
h(e_{1},e_{1}),h(e_{k},e_{k})\rangle]\nonumber\\
&\leq&\frac{n(n-1)}{2}(c+H^2)-\frac{3R}{2}+\sum_{k=2}^nRic(e_k)+(n-1)c.
\end{eqnarray}
 Therefore, if $$3R-2Ric^{[n-1]}_{\max}>n(n-1)(c+H^2),$$  then
$$\sum_{k=2}^{n}[2|h(e_{1},e_{k})|^{2}-\langle
h(e_{1},e_{1}),h(e_{k},e_{k})\rangle]<(n-1)c.$$
Hence, it follows from Theorem 2.1 that $H_{1}(M;\mathbb{Z})=H_{n-1}(M;\mathbb{Z})=0,$
 and  $\pi_1(M)=0$. This proves Lemma 2.2.\qed \\

 If $$Ric_M>\frac{n(n-1)}{n+2}(c+H^2),$$  then we can directly obtain
$$2Ric_{\min}+R>n(n-1)(c+H^2).$$ Therefore, we have the following corollary  from Lemma 2.2.\\\\
 \textbf{Corollary 2.3 (\cite{XG2}).} \emph{Let $M$ be an $n(\geq4)$-dimensional
 compact submanifold in $F^{n+p}(c)$ with  $c\geq0$.
 If the Ricci curvature of M satisfies
  $$Ric_M>\frac{n(n-1)}{n+2}(c+H^2),$$ then $H_{1}(M;\mathbb{Z})=H_{n-1}(M;\mathbb{Z})=0,$ and $\pi_1(M)=0$.}\\

 \section{Proof of the Main Theorem}
In this section, we will give the proof of our Main Theorem. we first need
to eliminate the homology groups
$H_{q}(M;\mathbb{Z})$ for $1<q<n-1$ under the s-th weak Ricci curvature and scalar curvature pinching condition.\\\\
\textbf{Lemma 3.1.} \emph{Let $M$ be an $n(\geq4)$-dimensional
compact submanifold in $F^{n+p}(c)$ with  $c\geq0$.
 Assume that one of the following  conditions holds: \\
  $(i)$ $(n-q)(q-1)R+(n-2q)Ric^{[q]}_{\min}
> n^2(n-q-1)(q-1)(c+H^2),$ \\
$(ii)$  $q(n-q-1)R-(n-2q)Ric^{[n-q]}_{\max}
> n^2(n-q-1)(q-1)(c+H^2),$\\
where q is some integer in $[2, \frac{n}{2}]$. Then $H_{q}(M;Z)=H_{n-q}(M;Z)=0$.}\\\\
\textbf{Proof.} Assume that $ 2\leq q\leq \frac{n}{2}$. Setting
$$T_{\alpha}:=\frac{trH_{\alpha}}{n},$$
  we have $\sum_{\alpha}T_{\alpha}^{2}=H^2.$  Then we get from (\ref{Ricci})
  \begin{equation}Ric(e_i)=(n-1)c+\sum_{\alpha}\Big[nT_{\alpha}h^{\alpha}_{ii}-(h^{\alpha}_{ii})^2
  -\sum_{j\neq i}(h^{\alpha}_{ij})^2\Big].\label{3.1}
  \end{equation}
Then we get\begin{eqnarray}
&&\sum_{k=q+1}^{n}\sum_{i=1}^{q}[2|h(e_{i},e_{k})|^{2}-\langle
h(e_{i},e_{i}),h(e_{k},e_{k})\rangle]\nonumber\\
&=&2\sum_{\alpha}\sum_{k=q+1}^{n}\sum_{i=1}^{q}(h^{\alpha}_{ik})^{2}-
\sum_{\alpha}\sum_{k=q+1}^{n}\sum_{i=1}^{q}h^{\alpha}_{ii}h^{\alpha}_{kk}\nonumber\\
&=&\sum_{\alpha}\Big[2\sum_{k=q+1}^{n}\sum_{i=1}^{q}(h^{\alpha}_{ik})^{2}-
a\Big(\sum_{i=1}^{q}h^{\alpha}_{ii}\Big)\Big(nT_{\alpha}-\sum_{i=1}^{q}h^{\alpha}_{ii}\Big)\nonumber\\
&&-(1-a)\Big(\sum_{k=q+1}^{n}h^{\alpha}_{kk}\Big)\Big(nT_{\alpha}-\sum_{k=q+1}^{n}h^{\alpha}_{kk}\Big)\Big]\label{3.2}
\end{eqnarray}
for some $a\in(0,1)$.
Then it follows from (\ref{3.1}) and the Cauchy-Schwarz inequality that
\begin{eqnarray}&&\sum_{k=q+1}^{n}\sum_{i=1}^{q}[2|h(e_{i},e_{k})|^{2}-\langle
h(e_{i},e_{i}),h(e_{k},e_{k})\rangle]\nonumber\\
&\leq&\sum_{\alpha}\Big[2\sum_{k=q+1}^{n}\sum_{i=1}^{q}(h^{\alpha}_{ik})^{2}-
anT_{\alpha}\sum_{i=1}^{q}h^{\alpha}_{ii}+tq\sum_{i=1}^{q}(h^{\alpha}_{ii})^2\nonumber\\
&&-(1-a)nT_{\alpha}\sum_{k=q+1}^{n}h^{\alpha}_{kk}+(1-a)(n-q)\sum_{k=q+1}^{n}(h^{\alpha}_{kk})^2\Big]\nonumber\\
&\leq&\sum_{\alpha}\Big[2\sum_{k=q+1}^{n}\sum_{i=1}^{q}(h^{\alpha}_{ik})^{2}-aq\sum_{k\neq i}\sum_{i=1}^{q}(h^{\alpha}_{ik})^{2}-
(1-a)(n-q)\sum_{i\neq k}\sum_{k=q+1}^{n}(h^{\alpha}_{ik})^{2}\Big]\nonumber\\
&&+ aq\sum_{i=1}^{q}[(n-1)c-Ric(e_i)]+(1-a)(n-q)\sum_{k=q+1}^{n}[(n-1)c-Ric(e_k)]
\nonumber\\
&&+an(q-1)\sum_{\alpha}\sum_{i=1}^{q}T_{\alpha}h^{\alpha}_{ii}
+(1-a)n(n-q-1)\sum_{\alpha}\sum_{k=q+1}^{n}T_{\alpha}h^{\alpha}_{kk}.\label{3.3}
\end{eqnarray}
Since $q\geq2$, and $n-q\geq2$, we get
\begin{eqnarray}&&\sum_{\alpha}\Big[aq\sum_{k\neq i}\sum_{i=1}^{q}(h^{\alpha}_{ik})^{2}+
(1-a)(n-q)\sum_{i\neq k}\sum_{k=q+1}^{n}(h^{\alpha}_{ik})^{2}\Big]\nonumber\\
&\geq&\sum_{\alpha}\Big[2a\sum_{k\neq i}\sum_{i=1}^{q}(h^{\alpha}_{ik})^{2}+
2(1-a)\sum_{i\neq k}\sum_{k=q+1}^{n}(h^{\alpha}_{ik})^{2}\Big]\nonumber\\
&\geq&\sum_{\alpha}2\sum_{k=q+1}^{n}\sum_{i=1}^{q}(h^{\alpha}_{ik})^{2}.\label{3.4}
\end{eqnarray}
Let $a=\frac{n-q-1}{n-2}$, then we have $a(q-1)=(1-a)(n-q-1).$
This together with (\ref{3.3}), (\ref{3.4}) implies
\begin{eqnarray}&&\sum_{k=q+1}^{n}\sum_{i=1}^{q}[2|h(e_{i},e_{k})|^{2}-\langle
h(e_{i},e_{i}),h(e_{k},e_{k})\rangle]-q(n-q)c\nonumber\\
&\leq& -(2a-1)\sum_{i=1}^{q}Ric(e_i)-(1-a)(n-q)R
\nonumber\\
&&+[aq^2+(1-a)(n-q)^2](n-1)c\nonumber\\
&&+an^2(q-1)\sum_{\alpha}T^2_{\alpha}-q(n-q)c\nonumber\\
&=&-\frac{n-2q}{n-2}\sum_{i=1}^{q}Ric(e_i)-\frac{(n-q)(q-1)}{n-2}R
\nonumber\\
&& +\frac{n^2(n-q-1)(q-1)}{n-2}(c+H^2)\label{3.5}
\end{eqnarray}
It follows from the condition (i) and (\ref{3.5}) that
\begin{eqnarray}&&\sum_{k=q+1}^{n}\sum_{i=1}^{q}[2|h(e_{i},e_{k})|^{2}-\langle
h(e_{i},e_{i}),h(e_{k},e_{k})\rangle]-q(n-q)c\nonumber\\
&\leq&-\frac{n-2q}{n-2}Ric^{[q]}_{\min}-\frac{(n-q)(q-1)}{n-2}R
\nonumber\\&& +\frac{n^2(n-q-1)(q-1)}{n-2}(c+H^2)\nonumber\\
&<&0.\label{3.6}\end{eqnarray}
We also get from (\ref{3.5}) and condition (ii) that
\begin{eqnarray}&&\sum_{k=q+1}^{n}\sum_{i=1}^{q}[2|h(e_{i},e_{k})|^{2}-\langle
h(e_{i},e_{i}),h(e_{k},e_{k})\rangle]-q(n-q)c\nonumber\\
&\leq&\frac{n-2q}{n-2}\sum_{k=q+1}^{n}Ric(e_k)-\frac{q(n-q-1)}{n-2}R\nonumber\\
&& +\frac{n^2(n-q-1)(q-1)}{n-2}(c+H^2)\nonumber\\
&\leq& \frac{n-2q}{n-2}Ric^{[n-q]}_{\max}-\frac{q(n-q-1)}{n-2}R
\nonumber\\&& +\frac{n^2(n-q-1)(q-1)}{n-2}(c+H^2)\nonumber\\
&<&0.\label{3.8}\end{eqnarray}
 Then the assertion follows from (\ref{3.6}), (\ref{3.8}) and Theorem 2.1.\qed \\\\
\textbf{ Proof of Theorem 1.1.} If $R>n(n-2)(c+H^2)$, and $Ric_M\geq \frac{n}{2}(c+H^2)$, then $$R+2Ric_{\min}>n(n-1)(c+H^2).$$
 It follows from Theorem 2.2 that $H_{1}(M;\mathbb{Z})=H_{n-1}(M;\mathbb{Z})=0,$ and $\pi_1(M)=0$.

For $q\in[2,\frac{n}{2}]$,  we consider the following three
    cases:

Case I. $q=\frac{n}{2}$. If $R>n(n-2)(c+H^2)$,  then
 $(n-q)(q-1)R+(n-2q)Ric^{[q]}_{\min}
> n^2(n-q-1)(q-1)(c+H^2)$. Hence we get from Theorem 3.1 that $H_{\frac{n}{2}}(M;Z)=0$.

Case II. $q<\frac{n}{2}$. If $R>n(n-2)(c+H^2)$ and $Ric^{[q]}_{\min}\geq n(q-1)(c+H^2)$, then
 $(n-q)(q-1)R+(n-2q)Ric^{[q]}_{\min}
> n^2(n-q-1)(q-1)(c+H^2),$ for some $q\in[2,\frac{n}{2})$. Therefore, we get from Theorem 3.1 that $H_{q}(M;Z)=H_{n-q}(M;Z)=0$ for some $q\in[2,\frac{n}{2})$.

Since   $Ric^{[l]}_{\min}\geq lRic_{\min}$, we have if
$n=even$, and $$Ric_M\geq (n-2-\frac{4}{n-2})(c+H^2),$$ then $Ric^{[q]}_{\min}\geq n(q-1)(c+H^2)$ for all $q\in[2,\frac{n}{2})$.
If $n=odd$, and $$Ric_M\geq (n-2-\frac{2}{n-1})(c+H^2),$$  then $Ric^{[q]}_{\min}\geq n(q-1)(c+H^2)$ for all $q\in[2,\frac{n}{2})$.

Moreover, we have $$\frac{n}{2}(c+H^2)\leq (n-2)(c+H^2)$$ for $n\geq4$,
 $$\frac{n}{2}(c+H^2)\leq (n-2-\frac{2}{n-1})(c+H^2)$$ for $n\geq 5$,
  and $$\frac{n}{2}(c+H^2)\leq (n-2-\frac{4}{n-2})(c+H^2)$$ for $n\geq 6$. Hence it follows from the assumption that $H_{q}(M;Z)=H_{n-q}(M;Z)=0$ for all $q\in[1,\frac{n}{2}]$ and $\pi_1(M)=0$. It follows from $H_{i}(M;Z)=0$ for $i=1,2,\cdots,n-1$ that $M$ is a homology sphere.
 Since $M$ is simply connected, $M$ is a homotopy sphere. This together with the generalized Poincar\'{e}
 conjecture implies that $M$ is a topological sphere. This completes the proof of Theorem 1.1. \qed\\

Moreover, when $M$ is oriented, we get the following sphere theorem.\\\\
\textbf{Theorem 3.2.} \emph{Let $M$ be an $n(\geq4)$-dimensional oriented compact submanifold in
$F^{n+p}(c)$ with  $c\geq0$. If  $R>n(n-2)(c+H^2)$, and $$Ric^{[2]}_{\min}\geq 2(n-2-\frac{2\sigma_n}{2n-\sigma_n})(c+H^2),$$
then $M$ is homeomorphic to a sphere. Here $\sigma_n=sgn(n-4)((-1)^n+3)$, and $sgn(\cdot)$ is the standard sign function.}\\\\
\textbf{Proof.}  We first prove that
$Ric_M>0$. We consider the following three cases:

$(i)$ If $n=4$, and if there exists a point $x_0\in M$ such that $Ric_{\min}(x_0)\leq 0$, then without loss of generality,
we assume that $Ric(e_1)\leq 0$ at point $x_0$. Hence it follows from the assumption that $Ric(e_i)\geq 2(n-2)(c+H^2)$ for $i=2,\cdots, n$.
Then we get \begin{eqnarray}R(x_0)&\geq& 2(n-2)(c+H^2)+(n-2)2(n-2)(c+H^2)\nonumber\\
&=&2(n-1)(n-2)(c+H^2)\nonumber\\
&=&n(n-1)(c+H^2).\end{eqnarray}
On the other hand, we know from the Gauss equation that $R\leq n(n-1)(c+H^2)$. Therefore $R(x_0)=n(n-1)(c+H^2)$, and $S(x_0)=nH^2$, i.e., $x_0$ is a totally umbilical point. Then we get from the Gauss equation that $Ric_M\equiv(n-1)(c+H^2)$ at that point, which contradicts with the assumption. Therefore, $Ric_M>0$.

$(ii)$ If $n$ is even, $n\geq 6$, and if there exists a point $x_1\in M$ such that $Ric_{\min}(x_1)\leq 0$, then without loss of generality,
we assume that $Ric(e_1)\leq 0$ at point $x_1$. Hence it follows from the assumption that $Ric(e_i)\geq 2(n-2-\frac{4}{n-2})(c+H^2)$ for $i=2,\cdots, n$.
Then we get \begin{eqnarray}R(x_1)&\geq& 2(n-2-\frac{4}{n-2})(c+H^2)+2(n-2)(n-2-\frac{4}{n-2})(c+H^2)\nonumber\\
&=&2(n-1)(n-2-\frac{4}{n-2})(c+H^2)\nonumber\\
&\geq&n(n-1)(c+H^2),\end{eqnarray}
and the equality holds only if $n=6$. This together with the Gauss equation implies that $R(x_1)=n(n-1)(c+H^2)$, and $n=6$.  But a similar argument as in Case (i) shows that $Ric_M\equiv(n-1)(c+H^2)$ at point $x_1$, and it contradicts with the assumption. Therefore, $Ric_M>0$.

(iii) If $n$ is odd, and if there exists a point $x_2\in M$ such that $Ric_{\min}(x_2)\leq 0$, then without loss of generality,
we assume that $Ric(e_1)\leq 0$ at point $x_2$. Hence it follows from the assumption that $Ric(e_i)\geq 2(n-2-\frac{2}{n-1})(c+H^2)$ for $i=2,\cdots, n$.
Then we get \begin{eqnarray}R(x_2)&\geq& 2(n-2-\frac{2}{n-1})(c+H^2)+2(n-2)(n-2-\frac{2}{n-1})(c+H^2)\nonumber\\
&=&2(n-1)(n-2-\frac{2}{n-1})(c+H^2)\nonumber\\
&\geq&n(n-1)(c+H^2),\end{eqnarray}
and the equality holds only if $n=5$. This together with the Gauss equation implies that $R(x_2)=n(n-1)(c+H^2)$, and $n=5$.  But a similar argument as in Case (i) shows that $Ric_M\equiv(n-1)(c+H^2)$ at point $x_2$, and it contradicts with the assumption. Therefore, $Ric_M>0$.

On the other hand, it follows from the proof of Theorem 1.1 and the assumption that $H_{q}(M;Z)=H_{n-q}(M;Z)=0$ for all $q\in[2,\frac{n}{2}]$.
This together with the universal
coefficient theorem implies that
 $H^{n-1}(M;\mathbb{Z})$ has no torsion, and hence neither does $H_{1}(M;\mathbb{Z})$
 by the Poincar$\acute{{\rm e}}$ duality. We also get from $Ric_M>0$ and
the Bonnet-Myers theorem that the fundamental group $\pi_{1}(M)$ is
finite. Then we have $H_{1}(M;\mathbb{Z})=0$. Hence,
 $H_{n-1}(M;\mathbb{Z})=0$. Denote by $\widetilde{M}$ the universal Riemannian covering of
$M$. We may consider $\widetilde{M}$ be a Riemannian submanifold of
$F^{n+p}(c)$, and hence $\widetilde{M}$ is a homology sphere.
Since $\widetilde{M}$ is simply connected, it is a topological
sphere, which together with a result of Sjerve \cite{Sjerve} implies
that $M$ is simply connected. Then $M$ is a homotopy sphere. This together with the generalized Poincar\'{e}
 conjecture implies that $M$ is a topological sphere. \qed \\

It's easy to get that if $Ric^{[2]}_{\min}>2(n-2)(c+H^2)$, then $R>n(n-2)(c+H^2)$. Hence we have the following corollary.\\\\
\textbf{Corollary 3.3.} \emph{Let $M$ be an $n(\geq4)$-dimensional oriented compact submanifold in
$F^{n+p}(c)$ with  $c\geq0$. If $Ric^{[2]}_{\min}>2(n-2)(c+H^2)$,
then $M$ is homeomorphic to a sphere.}\\

 For $n=4$, we also get the following differentiable theorem.\\\\
 \textbf{Theorem 3.4.} \emph{Let $M$ be a $4$-dimensional  compact oriented submanifold in
$F^{n+p}(c)$ with $c\geq0$. If  $R>8(c+H^2)$, and $Ric_M>0$, then $M$ is diffeomorphic to a sphere.} \\\\
 \textbf{Proof.} Since $R>8(c+H^2)$ is equal to $S<4c+8H^2$, it follows from a theorem due to Xu-Zhao \cite{XZ} that $M$ is diffeomorphic to a spherical space form.
 Moreover, we know from Lemma 3.1 that $H_{2}(M;\mathbb{Z})=0$. Using the same argument as in the proof of Theorem 3.2 implies that
  $H_{1}(M;Z)=H_{n-1}(M;Z)=0$, and $M$ is simply connected. Therefore, $M$ is diffeomorphic to a sphere.\qed \\\\
\textbf{Theorem 3.5.} \emph{Let $M$ be an $n(\geq5)$-dimensional
oriented compact submanifold in $F^{n+p}(c)$ with  $c\geq0$. If $n$ is odd,
and
$$Ric^{[2]}_{\min}>2(n-2-\epsilon_n)(c+H^2),$$  then M
is homeomorphic to a sphere. Here
$\epsilon_n =\frac{4}{n^3-2n^2-n-2}.$} \\\\
\textbf{Proof.} Since $R\geq\frac{n}{l}Ric^{[l]}_{\min}$, we get from Lemma 3.1 that
 if $$Ric^{[q]}_{\min}>\Big[q(n-1)-\frac{q^2(n-q)(n-2)}{q(n-2q)+n(q-1)(n-q)}\Big](c+H^2),$$
where $q$ is some integer in $[2, \frac{n}{2}]$, then $H_{q}(M;Z)=H_{n-q}(M;Z)=0$.
On the other hand,  we know that $\frac{Ric^{[l]}_{\min}}{l}\geq\frac{ Ric^{[2]}_{\min}}{2}$ for $l\geq 2$, and the function $$f(s)=\frac{s(n-s)}{s(n-2s)+n(s-1)(n-s)}$$ is strictly monotone decreasing for $2\leq s\leq  \frac{n}{2}$. Hence, if $$Ric^{[2]}_{\min}>\Big[2(n-1)-\frac{2s(n-s)(n-2)}{s(n-2s)+n(s-1)(n-s)}\Big](c+H^2),$$ for
some integer $s\in[2, \frac{n}{2}]$, then $H_{q}(M;Z)=H_{n-q}(M;Z)=0$ for all
 $q\in[2, s].$ Let $s=\frac{n-1}{2}$. Then we get from the assumption that $H_{q}(M;Z)=H_{n-q}(M;Z)=0$ for $q\in[2,\frac{n}{2}]$.
 Moreover, it follows from the proof of Theorem 3.2 that $Ric_M>0$. Using a same argument as in the proof of Theorem 3.2, we get that
  $H_{1}(M;Z)=H_{n-1}(M;Z)=0$, and $M$ is simply connected. Then $M$ is a homotopy sphere. This together with the generalized Poincar\'{e}
 conjecture implies that $M$ is a topological sphere. This proves Theorem 3.5.\qed\\\\
\textbf{Theorem 3.6.} \emph{Let $M$ be an $n(\geq5)$-dimensional compact submanifold in
$F^{n+p}(c)$ with  $c\geq0$.  Assume that the scalar curvature of M satisfies  $R>n(n-2)(c+H^2)$. Then we have\\
$(i)$ if $n=5, 6, 8$, and  $\frac{Ric^{[s]}_{\max}}{s}\leq[n-2+\frac{n-4}{2(n-1)}](c+H^2),$
then $M$ is homeomorphic to a sphere.\\
$(ii)$ if $n=7$,  $n\geq 9$,  and $\frac{Ric^{[s]}_{\max}}{s}\leq (n-2+\frac{2\tau_n}{2n+\tau_n})(c+H^2)$,
then $M$ is homeomorphic to a sphere. \\
Here $s$ is some integer in $[1, \frac{n+2}{2}]$, and $\tau_n=(-1)^n+3$.}\\\\
\textbf{Remark.} It follows from (ii) of  Example 1.3 that $\frac{Ric^{[s]}_{\max}}{s}=(n-2+\frac{4}{n+2})(c+H^2)$, for the submanifold $S^{m-1}\big(\sqrt{\frac{m-1}{n(c+H^2)}}\big)\times
S^{m+1}\big(\sqrt{\frac{m+1}{n(c+H^2)}}\big)$ with $n=2m$. It follows from (iii) of  Example 1.3 that  $\frac{Ric^{[s]}_{\max}}{s}=(n-2+\frac{2}{n+1})(c+H^2),$ for the submanifold $S^{m}\big(\sqrt{\frac{m}{n(c+H^2)}}\big)\times
S^{m+1}\big(\sqrt{\frac{m+1}{n(c+H^2)}}\big)$ with $n=2m+1$. \\\\
\textbf{ Proof.}  If $R>n(n-2)(c+H^2)$, and $Ric^{[n-1]}_{\max}\leq (n^2-2.5n)(c+H^2)$, then $$3R-2Ric^{[n-1]}_{\max}>n(n-1)(c+H^2).$$
 It follows from Theorem 2.2 that $H_{1}(M;\mathbb{Z})=H_{n-1}(M;\mathbb{Z})=0,$ and $\pi_1(M)=0$.

For $q\in[2,\frac{n}{2}]$,  we consider the following three
    cases:

Case I. $q=\frac{n}{2}$. If $R>n(n-2)(c+H^2)$,  then
  $q(n-q-1)R-(n-2q)Ric^{[n-q]}_{\max}
> n^2(n-q-1)(q-1)(c+H^2).$ Hence we get from Lemma 3.1 that $H_{\frac{n}{2}}(M;Z)=0$.

Case II. $q<\frac{n}{2}$. If $R>n(n-2)(c+H^2)$ and $Ric^{[n-q]}_{\max}\leq n(n-q-1)(c+H^2)$, then
$q(n-q-1)R-(n-2q)Ric^{[n-q]}_{\max}
> n^2(n-q-1)(q-1)(c+H^2),$ for some $q\in[2,\frac{n}{2})$. Therefore, we get from Lemma 3.1 that $H_{q}(M;Z)=H_{n-q}(M;Z)=0$ for some $q\in[2,\frac{n}{2})$.

Since   $\frac{Ric^{[l_1]}_{\max}}{l_1}\leq \frac{Ric^{[l_2]}_{\max}}{l_2}$ for $l_1\geq l_2$, we have if
$n=even$, and $$ Ric^{[\frac{n+2}{2}]}_{\max}\leq \frac{n^2}{2}(c+H^2),$$ then $Ric^{[n-q]}_{\max}\leq n(n-q-1)(c+H^2)$ for all $q\in[2,\frac{n}{2})$.
If $n=odd$, and $$Ric^{[\frac{n+1}{2}]}_{\max}\leq \frac{n(n-1)}{2}(c+H^2),$$  then $Ric^{[n-q]}_{\max}\leq n(n-q-1)(c+H^2)$ for all $q\in[2,\frac{n}{2})$.

Moreover, we have $$\frac{n(2n-5)}{2(n-1)}<\frac{n(n-1)}{(n+1)}$$ for $n=5$, and $$\frac{n(2n-5)}{2(n-1)}\geq\frac{n(n-1)}{(n+1)}$$ for $n\geq7$;
$$\frac{n(2n-5)}{2(n-1)}<\frac{n^2}{n+2}$$ for $n=6,8$, and $$\frac{n(2n-5)}{2(n-1)}\geq\frac{n^2}{n+2}$$ for $n\geq10$.
Hence it follows from the assumption that $H_{q}(M;Z)=H_{n-q}(M;Z)=0$ for all $q\in[1,\frac{n}{2}]$ and $\pi_1(M)=0$. It follows from $H_{i}(M;Z)=0$ for $i=1,2,\cdots,n-1$ that $M$ is a homology sphere.
 Since $M$ is simply connected, $M$ is a homotopy sphere. This together with the generalized Poincar\'{e}
 conjecture implies that $M$ is a topological sphere. This completes the proof of Theorem 3.6.\qed\\\\
 \textbf{Theorem 3.7.} \emph{Let $M$ be an $n$-dimensional compact oriented  submanifold in
$F^{n+p}(c)$ with  $c\geq0$. Assume that the scalar curvature of M satisfies  $R>n(n-2)(c+H^2)$. We have that \\
$(i)$ if $n=5$, and  $\frac{Ric^{[s]}_{\max}}{s}\leq(n-2+\frac{2}{n+1})(c+H^2),$
then $M$ is homeomorphic to a sphere.\\
$(ii)$ if $n=4,6,8$,  and $\frac{Ric^{[s]}_{\max}}{s}\leq(n-2+\frac{4}{n+2})(c+H^2)$, then $M$ is homeomorphic to a sphere.\\
Here $s$ is some integer in $[1, \frac{n+2}{2}]$.}\\\\
\textbf{Proof.} It follows from the assumption and the proof of Theorem 3.6 that  $H_{2}(M;Z)=H_{3}(M;Z)=0$ for $n=5$, and $$H_{q}(M;Z)=H_{n-q}(M;Z)=0$$ for $n=4$, $6$, $8$, where $q\in[2,\frac{n}{2}]$. Moreover, we have
\begin{eqnarray}
Ric_{\min}&\geq& R-(n-1)\frac{Ric^{[l]}_{\max}}{l}\nonumber\\
&>& n(n-2)(c+H^2)-(n-1)(n-2+\frac{4}{n+2})(c+H^2)\nonumber\\
&=&\frac{n(n-4)}{n+2}(c+H^2)\nonumber\\
&>&0.
\end{eqnarray}
Then a similar argument as in the proof of Theorem 3.2 shows that $$H_{1}(M;Z)=H_{n-1}(M;Z)=0,$$ and $M$ is simply connected. Therefore, $M$ is homeomorphic to a sphere. This proves Theorem 3.7.\qed\\\\
\textbf{Proof of Theorem 1.2.}  Choosing $s=1$, we get the conclusion from Theorems 3.6 and 3.7.\qed\\

Motivated by Theorem 1.1, Example 1.3, the convergence results of Brendle and Schoen for
Ricci flow \cite{Brendle, Brendle0, Brendle2}, and the differentiable sphere theorems for submanifolds with positive Ricci curvature proved in \cite{XG2, XT}, we propose the following conjecture.\\\\
\textbf{Conjecture 3.8.} \emph{Let $M$ be an $n(\geq 5)$-dimensional compact
submanifold in the space form $F^{n+p}(c)$ with $c+H^2>0$.
 If $n$ is odd, and  $$Ric_{M}>(n-2-\frac{2}{n-1})(c+H^2),$$ then
$M$ is diffeomorphic to $S^n$.}\\

To verify Conjecture 3.8,  we hope to prove the
following conjecture on the normalized Ricci flow.\\\\
\textbf{Conjecture 3.9.} \emph{Let $(M,g_0)$ be an
$n(\geq5)$-dimensional compact submanifold in an $(n+p)$-dimensional
space form $F^{n+p}(c)$ with $c+H^2>0$. If $n$ is odd, and the Ricci curvature of
$M$ satisfies
$$Ric_{M} >
(n-2-\frac{2}{n-1})(c+H^2),$$
then the normalized Ricci flow with initial metric
$g_0$
$$\frac{\partial}{\partial t}g(t) = -2Ric_{g(t)} +\frac{2}{n}
r_{g(t)}g(t),$$ exists for all time and converges to a constant
curvature metric as $t\rightarrow\infty$. Moreover, $M$ is
diffeomorphic to a spherical space form.}\\

 \medskip
\noindent Juan-Ru Gu\\
Department of Applied Mathematics, Zhejiang University of Technology,  Hangzhou
310023, China; \\
Center of Mathematical Sciences,  Zhejiang University,  Hangzhou 310027, China\\
E-mail: gujr@zju.edu.cn\\

\noindent Hong-Wei Xu\\
Center of Mathematical Sciences, Zhejiang University,  Hangzhou 310027, China\\
E-mail: xuhw@zju.edu.cn

\medskip


\begin{thebibliography}{bb}


\bibitem{Chern}S. S. Chern, M. do Carmo, S. Kobayashi, Minimal submanifolds of a sphere with second
fundamental form of constant length, in Functional Analysis and
Related Fields, Springer-Verlag, New York(1970).
\bibitem{Brendle}S. Brendle, A general convergence result for the Ricci flow in higher dimensions,
 Duke Math. J. {\bf145}(2008) 585-601.
\bibitem{Brendle0}S. Brendle, Ricci Flow and the Sphere Theorem, Graduate Studies in Mathematics, Vol.{\bf111}, Americam Mathematical Society, 2010.
\bibitem{Brendle2} S. Brendle, R. Schoen, Manifolds with $1/4$-pinched curvature are space forms,
 J. Amer. Math. Soc.  {\bf22}(2009)  287-307.
\bibitem{Ejiri}N. Ejiri,  Compact minimal submanifolds of a sphere with
positive Ricci curvature,  J. Math. Soc. Japan  {\bf31}(1979)
251-256.
\bibitem{Gauchman}H. Gauchman, Minimal submanifolds of a sphere with bounded second fumdamental form,  Trans. Amer. Math. Soc. {\bf292}(1986) 779-791.
\bibitem{GX} J. R. Gu, H. W. Xu, The sphere theorems for
manifolds with positive scalar curvature,  J. Differential
Geom. {\bf92}(2012) 507-545.
\bibitem{HV} T. Hasanis, T. Vlachos, Ricci curvatures and minimal submanifolds, Pacific J. Math. {\bf197}(2001) 13-24.
\bibitem{Lawson}B. Lawson, Local rigidity theorems for minimal hyperfaces,  Ann. of Math. {\bf 89}(1969) 187-197.
\bibitem{Lawson2}B. Lawson, J. Simons, On stable currents and their
application to global problems in real and complex geometry,
 Ann. of Math.  {\bf98}(1973), 427-450.
\bibitem{Leung} P. F. Leung, Minimal submanifolds in a sphere II,  Bull. Landon Math. Soc. {\bf23}(1991) 387-390.
\bibitem{Li}A. M. Li, J. M. Li, An intrinsic rigidity theorem for minimal
submanifolds in a sphere, Arch. Math.  {\bf58}(1992)
582-594.
\bibitem{Li2}H. Z. Li, Curvature pinching for odd-dimensional minimal submanifolds in a
sphere,  Publ. Inst. Math. $($Beograd$)$  {\bf53}(1993)
122-132.
\bibitem{Micallef}M. Micallef, J. D. Moore, Minimal two-spheres and the topology of manifolds
with positive curvature on totally isotropic two-planes, Ann.
of Math. {\bf127}(1988) 199-227.
\bibitem{QT} C. Qian, Z. Z. Tang, Isoparametric foliations, a problem of Eells-Lemaire and conjectures of Leung,  Proc. London Math. Soc. {\bf112}(2016)  979-1001.
\bibitem{ShYB} Y. B. Shen, Curvature pinching for three-dimensional minimal
submanifolds in a sphere, Proc. Amer. Math. Soc. {\bf115}(1992) 791-795.
\bibitem {Shiohama1} K. Shiohama, H. W. Xu, The topological sphere theorem for complete submanifolds,
 Compositio Math. {\bf107}(1997), 221-232.
\bibitem{Shiohama2} K. Shiohama, H. W. Xu, An
integral formula for Lipschitz-Killing curvature and the critical
points for height functions, J. Geom. Anal. {\bf21}(2011) 241-251.
\bibitem{Simons}J. Simons, Minimal varieties in Riemannian manifolds, Ann. of Math. {\bf 88}(1968) 62-105.
\bibitem{Sjerve} D. Sjerve, Homology spheres which are covered by spheres, J. London Math. Soc. {\bf6}(1973) 333-336.
\bibitem{Sun}Z. Q. Sun, Submanifolds with constant mean
curvature in spheres, Adv. Math. $($China$)$ {\bf 16}(1987)
91-96.
\bibitem{Vlachos}T. Vlachos, A sphere theorem for odd-dimensional submanifolds of spheres,  Proc. Amer. Math. Soc. {\bf130}(2002) 167-173.
\bibitem{Xia}C. Y. Xia, A sphere theorem for submanifolds in a manifold with pinched positive curvature,  Monatsh Math.  {\bf124}(1997)  365-368.
\bibitem{Xin}Y. L. Xin, Application of integral currents to
vanishing theorems, Scient. Sinica(A) {\bf27}(1984)  233-241.
\bibitem{XG2}H. W. Xu, J. R. Gu, Geometric, topological and differentiable rigidity of
submanifolds in space forms, Geom. Funct. Anal.
{\bf23}(2013) 1684-1703.
\bibitem{XHZ} H. W. Xu, F. Huang, E. T. Zhao, Geometric and differentiable rigidity of submanifolds in spheres, J. Math. Pures Appl. {\bf99}(2013) 330-342.
\bibitem{XG3}H. W. Xu, Y. Leng, J. R. Gu, Geometric and topological rigidity for compact submanifolds of odd dimension, Sci. China Math. {\bf57}(2014) 1525-1538.
\bibitem{XT}H. W. Xu, L. Tian,  A differentiable sphere theorem inspired by
rigidity of minimal submanifolds, Pacific J. Math. {\bf254}(2011) 499-510.
\bibitem{XZ}H. W. Xu, E. T. Zhao, Topological and
differentiable sphere theorems for complete submanifolds, Comm.
Anal. Geom. {\bf 17}(2009) 565-585.
\bibitem{Yau}S. T. Yau, Submanifolds with constant mean curvature I, II, Amer. J. Math. {\bf96, 97}(1974, 1975) 346-366, 76-100.\\


\end{thebibliography}
\end{document}